\documentclass[a4paper,12pt]{article}
\usepackage{amsfonts,amssymb,latexsym}
\usepackage[cp1251]{inputenc}
\begin{document}

\author{A.~N.~Panov}
\date{\it{ Russia, 443011, Samara, ul.Akad. Pavlova 1,
Samara State University,
Mathematical Department\\
apanov@list.ru}}
\title{Representations of quantum orders} \maketitle \footnote{The paper is supported
by the RFFI-grants 02-01-00017 and 03-01-00167.}

\begin{abstract}
We study finite dimensional algebras that appear as  fibers of
quantum orders over a given point of  variety of center. We present
the formula for the number of irreducible representations and check
it for it for the algebra of twisted polynomials, the quantum Weyl
algebra and the algebra of regular functions on quantum group.
\end{abstract}

\section{ Introduction and main statements}

Quantum algebras appear in the framework of mathematical physics.
From the algebraic point of view a quantum algebra  $R_q$ is  an
domain and a free ${\Bbb C}[q,q^{-1}]$-module. After specialization
module $q-\varepsilon$ we obtain the algebra $R_\varepsilon =
R_q\bmod(q-\varepsilon)$. As usual (see algebras A1-A4 below)
$R_\varepsilon$ is a domain and, if $\varepsilon$ is a root of
unity, then $R_\varepsilon$ is  finite over its center
$Z_\varepsilon$. We call $R_\varepsilon$ a quantum order (since it
becomes an order in the skew field
$R_\varepsilon\otimes{\rm{Fract}}(Z_\varepsilon)$). This algebra
defines the affine variety ${\frak X}= {\rm{Maxspec}} Z_\varepsilon$
that is singular in general.

The order $R_\varepsilon$ has one remarkable property: it admits the
quantum adjoint action. For  $a,u\in R_q$ we denote
$a_\varepsilon,u_\varepsilon:=a, u \bmod(q-\varepsilon)\in
R_\varepsilon$. If $u_\varepsilon$ lies in $Z_\varepsilon$, then the
formula
$${\cal D}_u (a_\varepsilon) = \frac{ua -au}{q-\varepsilon}\bmod(q-\varepsilon)
$$ defines the derivation  ${\cal D}_u: R_\varepsilon \to R_\varepsilon$, that is called
the quantum adjoint action of  $u$ ~\cite{DC-K, DCKP, DC-P, PanTG,
PanAA}. The center $Z_\varepsilon$ is a Poisson algebra with respect
to the bracket $\{u_\varepsilon,v_\varepsilon\}:={\cal
D}_u(v_\varepsilon)=-{\cal D}_v(u_\varepsilon)$. The variety ${\frak
X}= {\rm{Maxspec}} Z_\varepsilon$ is a Poisson variety. It splits
into  symplectic leaves ~\cite{BGPoisson}.

For  a point $\chi\in {\frak X}={\rm{Maxspec}} Z_\varepsilon$, we
consider the finite dimensional algebra $R_\chi = R_\varepsilon/
R_\varepsilon m_\chi$. Call it the fiber of $R_\varepsilon$. An
irreducible representation of $R_\varepsilon$ with central character
$\chi$ passes through $R_\varepsilon\to R_\chi$. Therefore, these
representations are in one to one correspondence with irreducible
representations of $R_\chi$

The goal of this paper is to characterize the fibers in terms of a
point $\chi$ of Poisson variety ${\frak X}$. Main Theorem will be
proved in the case when $R_q$ is one of the following algebras: A1)
Algebra of twisted polynomials, A2) Quantum Weyl Algebra, A3)
$U_q({\frak b})$ (this algebra is isomorphic to ${\Bbb C}_q[B]$ for
the Borel subroup $B$), A4) algebra ${\Bbb C}_q[G]$ of regular
functions on the quantum semisimple Lie group $G$. For definitions
see e.g.~\cite{G}

We introduce the notion of stabilizer for any point of commutative
associative Poisson ${\Bbb C}$-algebra ${\cal F}$. Recall that an
algebra ${\cal F}$  is a Poisson algebra if it admits
 a Poisson bracket that is  a
linear, skew-symmetric map $\{\cdot,\cdot\}: {\cal F}\otimes {\cal
F}\to {\cal F}$, subjecting to the Jacobi and the Leibniz identity
(i.e., $\{a,bc\}=\{a,b\}c+b\{a,c\}$ for all $a,b,c\in {\cal F}$).
 A Poisson algebra is a Lie algebra with respect to
the Poisson bracket. We call $P$ an ideal(resp. Poisson ideal) of
Poisson algebra ${\cal F}$ if $P$ is an ideal of commutative
associative algebra ${\cal F}$ (resp. $P$ is an ideal of ${\cal F}$
and $\{P, {\cal F}\}\subset P$). We identify a point $\chi\in{\frak
X}  = {\rm{Maxspec}} {\cal F}$ with the character $\chi: {\cal F}\to
{\Bbb C}$. We use the notation $m_\chi$ for the corresponding
maximal ideal in ${\cal F}$. The subalgebra
$$
G_\chi:=\{a\in {\cal F}:~\{a, {\cal F}\}\in m_\chi\}
$$
is a Poisson subalgebra of ${\cal F}$. The ideal  $m^2_{\chi}$ is
contained in $G_\chi$ and
is a Poisson ideal in $G_\chi$.\\
{\bf Definition 1.1}. The finite dimensional Lie ${\Bbb C}$-algebra
$${\frak g}_\chi:= G_\chi/m_\chi^2$$
is called the stabilizer of the  point $\chi\in {\frak X}$. If
${\cal F}$ is generated (as commutative associative  ${\Bbb
C}$-algebra) by $a_1,\ldots, a_N$, then ${\frak g}$ is the linear
span of $\overline{a_i}:= a_i-\chi(a_i)\bmod m^2_\chi$. The
definition of stabilizer in the case of smooth manifolds is given in
~\cite{KM}.

If ${\frak g}$ is a finite dimensional Lie algebra over ${\Bbb C}$
and ${\frak n}$ is the maximal nilpotent ideal (i.e., nilradical) of
${\frak g}$, then the Lie algebra ${\frak g}/{\frak n}$ is a
reductive algebra. We denote by ${\rm{rank}}{\frak g}$ the dimension
of maximal commutative subalgebra in ${\frak g}/{\frak n}$. If
${\frak g}$ is an algebraic solvable Lie algebra (i.e., the Lie
algebra of some algebraic solvable ${\Bbb C}$-group ${\frak G}$),
then ${\frak g}={\frak t}\oplus {\frak n}$ where ${\frak t}$ is a
toroidal Lie subalgebra of ${\frak g}$. In this case $\dim_{\Bbb C}
{\frak t}= {\rm{rank}}{\frak g}$.

Recall the definition of Ore extension. Given  an algebra $A$; an
automorphism $\tau:A\to A$ and a $\tau$-derivation $\delta:A\to A$
(i.e. $\delta(ab)= \delta(a)b+\tau(a)\delta(b)$ for all $a,b\in A$).
An algebra $R$ is an Ore extension of $A$ if $R$ is generated by $A$
and an indeterminate $x$ with the defining relations
$xa=\tau(a)x+\delta(a)$ for all $a\in A$ ~\cite{G, MC-R}.

Let we have a skew-symmetric integer matrix ${\Bbb S} =
(s_{ij})_{ij=1}^N$ and an indeterminate $q$. Put $q_{ij}=q^{s_{ij}}$
and form the matrix ${\Bbb Q}=(q_{ij})_{ij}^n$. By definition, an
algebra $R_q$ is a quantum solvable algebra over $C:={\Bbb
C}[q,q^{-1}]$ if it is generated by the elements $x_1, x_2,\ldots,
x_n, x_{n+1}^{\pm 1},\ldots, x_{N}^{\pm 1}$ and $C$ and such that
any  its subalgebra
$$R_{i} := <x_i,\ldots, x_n, x_{n+1}^{\pm 1},\ldots,
x_{N}^{\pm 1}, C >,\quad 1\le i\le n$$
  is an Ore extension $R_i=R_{i+1}[x_i;\tau_i,\delta_i])$
  with $\tau(x_j)=q_{ij}x_j$, ~$~ i+1\le j\le N$ and
  $x_ix_j=q_{ij}x_jx_i$, $~1\le i\le N,~n+1\le j\le N$ ~\cite{PanA,
PanAA}. All algebras A1-A4 are quantum solvable algebras (see
~\cite{PanA, PanAA}). More precisely,  the algebra ${\Bbb C}_q[G]$
becomes quantum solvable after some localization (see section 3). In
what follows we
suppose that all quantum solvable algebras obey the following conditions. \\
1) "$q$-skew condition":
$\tau_i\delta_i=q^{s_i}\delta_i\tau_i{~\rm{with~some}~} s_i\in{\Bbb
Z}$; assume that $s_i\ne 0$ for $\delta_i\ne 0$ and $s_i=0$ for
$\delta_i=0$.
Call the system of non-zero integers $\{s_i\} $ the system  of exponents of $R_q$;\\
2) All $\delta_i$ is locally nilpotent.\\
Notice that the algebras A1-A4 obey these two conditions.

 Let $\varepsilon$ be a primitive $l$th root of unity.
 For the algebras A1-A2 we call $l$ (and
$\varepsilon$)  admissible if
 if $l$ is relatively prime with all principal minors of ${\Bbb S}$
and with the system of exponents $s_1,\ldots,s_N$. For the algebras
A3-A4 $l$(and $\varepsilon$) is admissible if $l$ is odd and $l\ge
3$ in the case $G$ has ${\rm{G}}_2$ components.

If $R_q$ is one of the algebras A1-A4 and $l$ is an admissible root
of unity, then the elements $a_i:= x_{i,\varepsilon}^l~ 1\le i\le N$
(here $x_{i,\varepsilon}=x_i\bmod(q-\varepsilon)$) lie in the center
$Z_\varepsilon$  of $R_\varepsilon$~\cite[Lemma 2.19]{PanAA}. The
central subalgebra $Z_0$, generated by $a_i,~ 1\le i\le n$ and
$a_i^{-1},~n+1\le i \le N$, is isomorphic to ${\Bbb C}[a_1,\ldots,
a_n, a_{n+1}^{\pm 1}, \ldots a_{N}^{\pm 1}]$. The subalgebra $Z_0$
is called the $l$-center of $R_\varepsilon$. For all algebras A1-A4
the subalgebra $Z_0$ is a Poisson subalgebra of $Z_\varepsilon$. The
embedding $Z_0\subset Z_\varepsilon$ defines the projection
$\phi:{\frak X}\to{\frak X}_0$ where ${\frak X}_0={\rm{Maxspec}}
Z_0$.

The goal of this paper is to prove the following statement for the
algebras
A1-A4.\\
{\bf Main Theorem.} Let $R$ be one of the algebras A1-A4. Suppose
that $l$ is an admissible. Let $\phi:\chi\mapsto\chi_0$ and ${\frak
g}_\chi$ (resp. ${\frak g}_{\chi_0})$ are stabilizer of $\chi$(resp.
$\chi_0$).
Then\\
1)  ${\frak g}_\chi$ and ${\frak g}_{\chi_0}$ are algebraic solvable
Lie algebras with decompositions ${\frak g}_\chi={\frak t}\oplus
{\frak n}$ and
${\frak g}_{\chi_0} = {\frak t}_0\oplus {\frak n}_0$.\\
2) the subalgebra $G_{\chi_0}$ (of the algebra $Z_0$) is contained
in $G_\chi$; the embedding $G_{\chi_0}\subset G_{\chi}$ is extended
to the homomorphism $\psi_p:{\frak g}_{\chi_0}\to {\frak g}_\chi $
such that its restriction on ${\frak t}_0$ is an isomorphism
$\psi_p\vert_{{\frak t}_0}:{\frak t}_{\chi_0}\to {\frak t}_\chi $;
${\rm{rank}} {\frak g}_\chi = {\rm{rank}}{\frak g}_{\chi_0}$ ;\\
3) the number $\vert\rm{Irr}R_\chi\vert$ of irreducible
representations of $R_\varepsilon$ with central character $\chi$  is
equal to $l^{{\rm{rank}} {\frak g}_{\chi}}.$

In the  section 2 we prove  Main Theorem in the partial case
(Proposition 2.2) and show that  Main Theorem is true for the
quantum solvable algebras with "admissible stratification"
(Proposition 2.5). In the next section 3 we prove the existence of
admissible stratifications for the algebras A1-A4. This will
conclude the proof of  Main Theorem for the algebras A1-A4 (see
Propositions 3.1-3.4).

\section{Standard ideals}

Let ${\Bbb S} = (s_{ij})_{ij=1}^N$, $C:={\Bbb C}[q,q^{-1}]$ be as
above. We denote by $A_{\Bbb Q}$ the algebra of twisted polynomials
(see section 3).

  Consider a quantum solvable algebra ${\frak R}'$ over $C$, generated by
$$x_1, \ldots, x_m, x_{m+1}, \ldots, x_N.$$
Suppose that the elements $x_1,\ldots,x_m$ are $q$-commute, i.e.
$x_ix_j=q_{ij}x_jx_i$, $1\le i,j\le m$. The multiplicatively closed
subset ${\frak S}$, generated  by   $x_1,\ldots,x_m$,  is  a
denominator subset ~\cite[Lemma 2.1]{C}.  Denote ${\frak R}:={\frak
R}'{\frak S}^{-1}$.
 There exist the elements
  $\tilde{x}_{m+1},\ldots,\tilde{x}_{N}$ in the localization of ${\frak R} S_*^{-1}$
  (here $S_*$ is  finitely generated, by some $\{q^i-1\}$, denominator subset
   in $C$) such that $x_i\tilde{x}_{j}=q_{ij}\tilde{x}_{j}x_i$, $1\le i\le
m$, $m+1\le j\le N$. The subalgebra, generated by
 $\tilde{x}_{m+1},\ldots,\tilde{x}_{N}$ coincides with
 ${\frak R}_{m+1} = < x_{m+1},\ldots, x_N >$~\cite[Prop.2.1-2.3]{C}.
 Suppose that the ideal ${\frak I}'$ ( of ${\frak R} S_*^{-1}$), generated by
$\tilde{x}_j$, $m+1\le j\le N$, has zero intersection with $C$.
Denote ${\frak I}:={\frak I}'\cap {\frak R}$. Call $({\frak
R},{\frak I})$ a standard pair and ${\frak I}$ a standard ideal in
${\frak R}$.
 Notice that the subalgebra ${\frak B}$, generated over
$C$ by $x_i, 1\le i\le m$, is an algebra of twisted Laurent
polynomials and ${\frak B}={\frak R}/{\frak I}$.

Let $\varepsilon$ be a primitive $l$th root of unity and $l$ is
relatively prime with  all principal minors of ${\Bbb S}$.
 We denote by
${\frak Z}_\varepsilon$ the center of ${\frak R}_\varepsilon:={\frak
R}\bmod(q-\varepsilon)$ and ${\frak X}:= {\rm{Maxspec}}{\frak
Z}_\varepsilon$. The elements $a_i:=x_{i,\varepsilon}^l, 1\le i\le
m$ lie in ${\frak Z}_\varepsilon$. Let ${\frak Z}_0$ be some
subalgebra of ${\frak Z}_\varepsilon$ such that ${\frak Z}_0\cap
{\frak B}_\varepsilon $ is generated by $a_i:=x_{i,\varepsilon}^l,
1\le i\le m$ and ${\frak R}_\varepsilon$ is finite over ${\frak
Z}_0$. The center ${\frak Z}_\varepsilon$ is finite over ${\frak
Z}_0$.
 Denote ${\frak X}_0:= {\rm{Maxspec}}{\frak Z}_0$, ${\frak I}_\varepsilon=({\frak I}+{\frak R}(q-\varepsilon))\bmod(q-\varepsilon)$,
 $\phi: {\frak X} \to {\frak X}_0$,
${\frak i}:={\frak I}_\varepsilon\cap {\frak Z}_\varepsilon$,
${\frak i}_0:={\frak I}_\varepsilon\cap {\frak Z}_0$.

The skew field ${\rm{Fract}}({\frak R})$ is isomorphic to the skew
field ${\rm{Fract}}(A_{\Bbb Q})$ of the algebra $A_{\Bbb
Q}$(=${\rm{gr}}({\frak R})$) of twisted polynomials (see section
3)~\cite{PanA, PanAA, C}.
 We are going to prove Main Theorem for the case $\chi$ (resp.$\chi_0$)
  is a point of  ${\frak X}$ (resp. ${\frak X}_0$) annihilated by
${\frak i}$ (resp.${\frak i}_0$).

 The algebra ${\frak B}$ has a new system of generators $h_i, g_i, ~1\le
i\le k$ and $z_j, ~1\le j\le p$, $2k+p=m$ (that consists of
monomials of $x^{\pm 1}_i$, $1\le i\le m$)  such that
$h_ig_i=q^{d_i'}g_ih_i$ and $\{z_j\}$ generate the center of ${\frak
B}$. By assumption, $l$ is relatively prime with $d'_i$.
 The
intersection of the center ${\frak Z}$ of ${\frak R}$ with ${\frak
B}$ is generated by some monomials $\{z^a:=z_1^{\alpha_1}\cdots
z_p^{\alpha_p},~ \alpha_j\in {\Bbb Z}\}$.
 Choosing the compatible basis, we may consider
that  the intersection ${\frak Z}\cap {\frak B}$ is generated by
$z_{t+1}^{n_{t+1}},\ldots, z_p^{n_p}$ for some
$n_{t+1},\ldots,n_p\in{\Bbb N}$. Since the  field
$\rm{Center}({\rm{Fract}}((A_{\Bbb Q}))$ is algebraically closed in
${\rm{Fract}}(A_{\Bbb Q})$, then if an element $z^d$ also lies in
the center of ${\rm{Fract}}(A_{\Bbb Q})$(that is isomorphic to
${\rm{Fract}}({\frak R})$), then $z$
lies in the center. This verifies that $n_{t+1}=\cdots=n_p=1$.\\
{\bf Lemma 2.1}. \\
1) The intersection ${\frak Z}_\varepsilon\cap {\frak
B}_\varepsilon$ is generated by
$$h^l_{i,\varepsilon},\quad g^l_{i,\varepsilon},{~{\rm{where}}~} 1\le i\le k {~\rm{and}~ }
z^l_{1,\varepsilon},\ldots, z^l_{t,\varepsilon},~
z_{t+1,\varepsilon},\ldots, z_{m,\varepsilon}.$$ 2)  The
intersection ${\frak Z}_0\cap {\frak B}_\varepsilon$ is generated by
$$h^l_{i,\varepsilon},\quad g^l_{i,\varepsilon}{~{\rm{where}}~ } 1\le i\le k   {~\rm{and}~}
z^l_{1,\varepsilon},\ldots, z^l_{t,\varepsilon},~
z^l_{t+1,\varepsilon},\ldots, z^l_{m,\varepsilon}.$$
{\bf Proof}.
The statement 2) is trivial. To prove 1) it suffices to show that
the monomial $z_{1,\varepsilon}^{\alpha_1}\ldots
z_{t,\varepsilon}^{\alpha_t}$ lies in ${\frak Z}_\varepsilon$
whenever $l$ divides all $\alpha_i$.

There exists the system of generators $\tilde{x}_{k_1},\ldots,
\tilde{x}_{k_t}$ in $A_{\Bbb Q}$
 such that
$$ z_j\tilde{x}_{k_j} = q^{\nu_{i,k_j}}\tilde{x}_{k_j}z_j \quad \rm{and}
~F:= \det(\nu_{i,k_j})_{i,j=1}^t\ne 0.$$ Connect
$\tilde{x}_{k_1},\ldots, \tilde{x}_{k_t}$ to the system $\{x_i,1\le
i\le m\}$. Denote by ${\Bbb S}''$ the corresponding
$(m+t)\times(m+t)$-submatrix of ${\Bbb S}$. The rank of ${\Bbb S}''$
is equal to $2k+2t$ and the greatest common divisor $D''$ of all its
$(2k+2t)\times(2k+2t)$-minors is equal to $ (d'_1)^2\ldots (d'_k)^2
F^2$. Since $l$ is admissible, $l$ relatively prime with $D''$.
Therefore, $GCD(l,F)=1$. There exist $v_i\in {\rm{Fract}}({\frak
R})$  such that
$$z_iv_j = q^{p_i\delta_{ij}}v_jz_i,~{\rm{and}}~{\rm{GCD}}(l,p_i)=1$$
for all $1\le i,j \le t$. This implies that, if
$z_{1,\varepsilon}^{\alpha_1}\ldots z_{t,\varepsilon}^{\alpha_t}$
lies in ${\frak Z}_\varepsilon$, then  $l$ divides all
$\alpha_i$.$\Box$

{\bf Proposition 2.2}. Let ${\frak R}$, ${\frak B}$, ${\frak I}$,
$\varepsilon$ be as above. Suppose that ${\frak Z}_0$ is a Poisson
subalgebra in ${\frak Z}_\varepsilon$. Let $\chi\in{\frak X}$ and
$\chi_0=\phi(\chi)\in{\frak X}_0$. Suppose that $\chi$ (resp.
$\chi_0$) is annihilated by the ideal ${\frak i}$ (resp.${\frak
i}_0$) and $\chi(a_i)\ne 0$, $1\le i\le m$.
Then \\
1) the number of irreducible representations of ${\frak
R}_\varepsilon$ with central character $\chi$ is equal to
$l^t$;\\
2) the ideal  ${\frak i}$ (resp. ${\frak i}_0$)  is a Poisson ideal
in $G_\chi$ (resp. $G_{\chi_0}$). Denote by ${\frak n}'$ (resp.
${\frak n}_0'$) the image of ${\frak i}$ (resp. ${\frak i}_0$) in
${\frak g}_\chi$ (resp. ${\frak g}_{\chi_0}$);\\
3) the ideal  ${\frak n}'$ (resp. ${\frak n}'_0$) is a nilpotent
ideal in ${\frak g}_\chi$ (resp. ${\frak g}_{\chi_0}$). Then Main
Theorem is true for ${\frak R}$ and $\chi$. In particular,  ${\frak
g}_\chi$ (resp. ${\frak g}_{\chi_0}$)
is an  algebraic solvable Lie algebra.\\
{\bf Proof}. First, notice that ${\frak I}$ lies in the radical of
${\frak R}_\varepsilon{\frak i}_0$ (apply~\cite[Lemma 5.1]{PanAA}).
Kernel of any irreducible representation $\pi$  with $l$-central
character $\chi_0$ contains ${\frak I}$. Any irreducible
representation  with $l$-central character $\chi_0\in {\frak X}_0$
is uniquely determined by its kernel, generated by
$${\frak i},\quad
h_{i,\varepsilon}^l -\chi(h_{i,\varepsilon}^l),~g_{i,\varepsilon}^l
- \chi(g_{i,\varepsilon}^l)~{\rm{for}}~1\le i\le k,$$
$$z_{j,\varepsilon} -
\chi(z_{j,\varepsilon}) ~{\rm{for}}~ 1\le j\le t,\quad
z_{j,\varepsilon} -\chi(z_{j,\varepsilon})~{\rm{for}}~ t+1\le j\le
p.$$ The
 number of irreducible representations
with central character $\chi$ is equal to $l^t$. This proves 1).

To calculate subalgebras ${\frak g}_\chi$ and ${\frak g}_{\chi_0}$
we find generators of the subalgebras $G_\chi$ and $G_{\chi_0}$ of
${\frak R}_\varepsilon$:

$$G_\chi = < z^l_{j,\varepsilon} - \chi(z^l_{j,\varepsilon}), 1\le j\le t;~
z_{j,\varepsilon} - \chi(z_{j,\varepsilon}), t+1\le j\le p;~~ {\frak
i}
>,$$
$$G_{\chi_0} = < z^l_{j,\varepsilon} - \chi(z^l_{j,\varepsilon}), 1\le j\le p;~
{\frak i}_0 >.$$ We see $G_{\chi_0}\subset G_\chi$. This defines the
homomorphism $\psi:{\frak g}_{\chi_0}\to{\frak g}_\chi$.

Since   ${\frak i}$ is an intersection of ${\frak R}$-ideal ${\frak
I}$ with ${\frak Z}_\varepsilon$, then ${\frak i}$ is a Poisson
ideal of ${\frak Z}_\varepsilon$ ~\cite[Lemma 3.12]{PanTG}. Similar
for ${\frak i}_0$. This implies 2). Denote

$$e_i = z^l_{i\varepsilon} - \chi(z^l_{i\varepsilon}) \bmod m_\chi^2,\quad 1\le i\le t
~\rm{and}~{\frak t}=\rm{span}\{e_i; ~1\le i\le t\}.$$ The Lie
algebra ${\frak g}_\chi$ is a  sum (as a linear space) of ${\frak
t}$ and ${\frak n}$
 spanned modulo $m^2_\chi$ by
$z_{j,\varepsilon} - \chi(z_{j,\varepsilon})$, $t+1\le j\le p $ and
${\frak n}'$.
 Similarly, ${\frak g}_{\chi_0}$ is a sum (as a linear space) of
${\frak t}_0={\frak t}$ and ${\frak n}_0$, spanned modulo
$m^2_{\chi_0}$ by  $z^l_{j,\varepsilon} -
\chi(z^l_{j,\varepsilon})$, $t+1\le j\le p$ and
 ${\frak n}'_0$.

Let us prove that  ${\frak n}$  is a nilpotent ideal in ${\frak
g}_\chi$. Similar for ${\frak n}_0$. Any element of ${\frak R}
S_*^{-1}$ (is a sum of monomials of
$$x_1^{n_1}\cdots
x_m^{n_m}\tilde{x}_{m+1}^{n_{m+1}}\cdots\tilde{x}_N^{n_N}$$

We define the degree ${\rm{deg}}(a):=(n_{m+1},\ldots,n_N)$ of any
monomial $a$. For any two monomials $a,b$ there exists $s\in{\Bbb
Z}$ such that $ab-q^sba$ is a sum of monomials of lower degree with
respect to the lexicographical ordering. For any $A, B\in
Z_\varepsilon$ we have $\{A,B\}={\rm{const}}AB+\{the
~lower~terms\}$. This verifies that ${\frak n}$ is a nilpotent
ideal.

Let us prove that the Lie subalgebra ${\frak t}$ is diagonalizable.
The elements $x_1,\ldots,x_m$ are FA-elements in ${\frak
R}$~\cite{PanA,PanAA}. That is for any $1\le i\le m$ and $a\in
{\frak R}$ there exists a polynomial $f(t)$ (with roots in
$\{q^s\}_{s\in{\Bbb Z}}$) such that $f({\rm{Ad}}_{x_i}(a)=0$. The
adjoint action ${\rm{Ad}}_{x_i}$ is diagonalizable ~\cite{PanA}. One
can choose $f(t)$ with different roots
$q^{\gamma_1},\ldots,q^{\gamma_k}$. The  derivation ${{\cal
D}}'_i:=x_\varepsilon^{-l}{\cal D}_{x^l}:{\frak
R}_\varepsilon\to{\frak R}_\varepsilon$ for $x=x_i$ obey $f_1({{\cal
D}}'_i)(a_\varepsilon)=0$ where $f_1(t)$ is a polynomial with
different roots $c\gamma_1,\ldots,c\gamma_k$,
$c=l\varepsilon^{l-1}$. This imply that ${{\cal D}}'_i$ is
diagonalizable. The same is true for
$z_i^l$. Finally, ${\rm{ad}}_{e_i}$  are simultaneously diagonalizable.$\Box$ \\
{\bf Definition 2.3}. Let $R$ be domain with unit. Consider the set
of pairs $\{({\cal P}_\mu,S_\mu)\}$ where $S_\mu$ is a denominator
subset $R$ and ${\cal P}_\mu$ is a prime ideal  in $R$ (i.e ${\cal
P}_\mu\in{\rm{Spec}}(R)$) with empty intersection with $S_\mu$.
 We call $\{({\cal P}_\mu,S_\mu)\}$ a stratification of ${\rm{Spec}}(R)$
 if for any $I\in{\rm{Spec}}(R)$ there exists a unique $\mu$ such that
 $I\supset {\cal P}_\mu$ and $I\cap S_\mu=\emptyset$.
 If $R$ is a free $C$-module over commutative ring $C$,
 we assume, in addition,  that $I$ and any ${\cal P}_\mu$ have
 zero intersection with $C$. \\
 {\bf Definition 2.4}. Let $R_q$ be a quantum solvable algebras over
  $C:={\Bbb C}[q,q^{-1}]$ and $\{({\cal P}_\mu,S_\mu)\}$
  be a stratification of $R_q$.
We call   $\{({\cal P}_\mu,S_\mu)\}$  an admissible stratification
if
\\
1) for any $\mu$ there exists isomorphism
$\theta_\mu:R_qS_\mu^{-1}\to {\frak R}_\mu$ such that ${\frak
R}_\mu$ and
${\frak I}_\mu:=\theta_\mu({\cal P}_\mu)$ form a standard pair;\\
2)  the stratification $\{({\cal P}_\mu,S_\mu)\}$
 admits specialization modulo $q-\varepsilon$ (i.e. $\{({\cal P}_{\mu,\varepsilon},S_{\mu,\varepsilon})\}$ is a
 stratification of ${\frak R}_\varepsilon$).\\
3) $S_{\mu,\varepsilon}:= S_\mu\bmod(q-\varepsilon)\subset Z_0$ and
$\theta_\mu(S_\mu)$ is generated by $x_1^l,\ldots, x_m^l$.\\
 {\bf
Proposition 2.5}. Let $R_q$ be a quantum solvable algebra and $l$ be
admissible for $R_q$. Suppose that $x_{1,\varepsilon}^l,\ldots,
x_{N,\varepsilon}^l$ lie in the center of $R_\varepsilon$ and
generate a Poisson central  subalgebra (denote $Z_0$). Suppose that
$R_q$ has an admissible stratification $\{({\cal P}_\mu,S_\mu)\}$.
Then  Main Theorem is true for $\chi$.\\
{\bf Proof.} Let $\chi\in{\frak X}$.
 Choose $\mu$ such that $\chi(S_{\mu,\varepsilon})\ne 0$ and
$\chi$ is annihilated by ${\frak i}_\mu$. Apply Proposition 2.2. $\Box$\\

\section{Existence of admissible stratification}
To prove Main Theorem we present an admissible stratification for quantum algebras
A1-A4.\\
{\bf A1) The algebra of twisted polynomials}. Let the matrices
${\Bbb Q}$ and ${\Bbb S}$ be as in the above. The algebra  $R =
A_{\Bbb S}$ of twisted  polynomials
 is generated by $x_1,\ldots x_n, x_{n+1}^{\pm 1},\ldots, x_N^{\pm 1}$
 subject to the relations $x_ix_j=q_{ij}x_jx_i$.

 Choose some subset $T\subset \Lambda = \{1,2,\ldots,n\}$.
Consider the ideal ${\cal P}_T$ generated by $\{x_i:~i\in T\}$ and
the denominator subset  $S_T$ generated by
$\{x^l_i:~i\notin T\}$. \\
{\bf Proposition 3.1.}
Main Theorem is true for the Algebra of twisted polynomials.\\
{\bf Proof}. The set of pairs $\{({\cal P}_T, S_T)\}$ is an
admissible stratification. By direct calculations,
$\{a_i,a_j\}=cs_{ij}a_ia_j$ where $a_i=x_{i,\varepsilon}^l$ and
$c=l\varepsilon^{l-1}$.  Apply Proposition 2.5. $\Box$

{\bf A2) Quantum Weyl algebra}. Let ${\Bbb S} = (s_{ij})_{ij}^n$ be
skew-symmetric integer matrix and $q$ be indeterminate. As above we
put $q_{ij}=q^{s_{ij}}$ and form the matrix ${\Bbb Q}
=(q_{ij})_{ij=1}^n)$. Given non-zero integers $s_1,\ldots, s_n$
define
$$q_1 = q^{s_1},\ldots,q_n = q^{s_n}.$$
We consider two new matrices. The first matrix ${\Bbb P}
=((p_{ij})_{ij=1}^n)$ with entries subject
 $p_{ii}=p_{ij}p_{ji}=1$ and such that
 $$p_{ij}=q_iq_{ij}, ~\rm{for}~ i<j .$$
 The second one
${\Bbb R} = (r_{ij})_{ij=1}^n)$ has entries
$$
r_{ij} = \left\{ \begin{array}{l}  q_{ji},~\rm{if}~ i<j,\\
q_i,~\rm{if}~ i=j,\\
p_{ji}= q_j q_{ji},~\rm{if}~ i>j.\\
\end{array}\right.$$
Form skew-symmetric integer matrix ${\Bbb T} = (t_{ij})_{ij=1}^n$
such that $p_{ij} =q^{t_{ij}}$ and  integer matrix ${{\Bbb U}} =
(u_{ij})_{ij=1}^n$ such that $r_{ij} = q^ {u_{ij}}$. Form matrices

$$
{\Bbb Q}^* = \left(\begin{array}{cc}
{\Bbb Q} & -{\Bbb R}\\
{\Bbb R}& {\Bbb P}
\end{array} \right),\quad
{\Bbb S}^* = \left(\begin{array}{cc}
{\Bbb S} & -{{\Bbb U}}\\
{{\Bbb U}}& {\Bbb T}
\end{array} \right).
$$
{\bf Definition 3.2}. The Quantum Weyl algebra $W$ is generated by
$y_1,\ldots, y_n, x_1,\ldots, x_n$ with the following relations
$y_iy_j =q_{ij}y_jy_i$, $x_ix_j = p_{ij}x_jx_i$, $
x_iy_j=r_{ij}y_jy_i$ for $i\ne j$ and
$$x_iy_i = q_iy_ix_i  + \sum_{k<i} (q_k-1)y_kx_k + 1.\eqno (3.1)$$

 The algebra $W$ is an quantum solvable algebra over
${\Bbb C}[q,q^{-1}]$ with the system of exponents $s_1,\ldots,s_n$.
Denote $h_k = y_kx_k$. The relations imply
$$\left(\begin{array}{l}
x_i\\y_i
\end{array}\right) h_k =
h_k \left(\begin{array}{l}
q_ix_i\\q_i^{-1}y_i
\end{array}\right)
\quad ~{\rm{for}}~ i<k, $$

$$\left(\begin{array}{l}
x_i\\y_i
\end{array}\right) h_k =
h_k \left(\begin{array}{l}
x_i\\y_i
\end{array}\right)
\quad ~{\rm{for}}~ i>k. $$
For any $1\le i\le n$ we denote $w_i = 1 + \sum_{k\le i}(q_k-1)y_kx_k$.

One can rewrite (3.1) as follows
$x_iy_i = q_iy_ix_i +w_{i-1}$.
The variables $x_i, y_i, w_i$  obey the relations
$$ y_iw_j = \left\{
\begin{array}{l}
q_i^{-1}w_jy_i,  ~{\rm{for}}~ i\le j,\\
w_jy_i,  ~{\rm{for}}~ i > j,
\end{array}
\right., \quad
x_iw_j = \left\{
\begin{array}{l}
q_i w_jx_i,  ~{\rm{for}}~ i\le j,\\
w_jx_i,  ~{\rm{for}}~ i > j,
\end{array}
\right..
$$

 By definition, $\varepsilon$ is an admissible $l$th root of unity if $l$ is
relatively prime with all principal minors of $S^*$ and with
$s_1,\ldots,s_n$. The elements $a_i:= x_{i,\varepsilon}^l, b_i:=
y_{i,\varepsilon}^l$ lie in the center of $Z_\varepsilon$ of
$W_{\varepsilon}$ and generate the central subalgebra $Z_0$.

Denote  $f_i= w_{i,\varepsilon}^l$, $1\le i \le n$. Similarly to
~\cite{JZ}, one can prove that
 there exists the chain of non-zero complex
 numbers $\gamma_1,\ldots, \gamma_{n-1}$ such that
 $$
 f_i = 1 + \sum_{k<i}\gamma_ka_kb_k.$$
This implies that $f_i\in Z_0$.
 By direct calculations,
 $\{a_i, a_j\} = \gamma t_{ij} a_ia_j$,
 $ \{b_i, b_j\} = \gamma s_{ij} b_ib_j$,
 $ \{a_i, b_j\} = \gamma u_{ij} a_ib_j$,
 $\{a_i, b_i\} = \gamma s_{i} a_i b_i + f_{i-1}$
 where $\gamma =l^2\varepsilon^{-1}$.
We see that $Z_0$ is a Poisson subalgebra in $Z_\varepsilon$.

Denote $\Lambda=\{1,\ldots,n\}$.
 We shall call a triple $T = (T_1,T_2,T_3)$ of subsets
 of $\Lambda$  an admissible triple if
 $T_1\subseteq T_2\subseteq T_3$ and the following property holds:
 if  $i\in T_2$ then  $i$ and $i-1\in T_3$.

 Consider the ideal ${\cal P}_T$ of $W$ generated by
 $x_i$, $y_j$ and $w_k$  with $i\in T_1$, $j\in T_2 $ and $k\in T_3$.
Form the denominator subset $S_T$ generated by  the following
$q$-commuting elements $\{x^l_i, i\in T_2 - T_1\}$, $\{y^l_i, i\in
\Lambda - T_2\}$ and $\{w^l_i, i\in \Lambda - T_3\}$. The subset
$S_T$ has empty intersection with ${\cal P}_T$. The set of pairs
$\{({\cal P}_T,S_T)\}$ is a stratification of $W$ (see
~\cite{OH,Horton}).

{\bf Proposition 3.3.}
Main Theorem is true for the Quantum Weyl algebra.\\
{\bf Proof}. The set of pairs $({\cal P}_T, S_T)$ is an admissible
stratification. $\Box$

 {\bf A3-A4. Cases of algebras
$U_q({\frak b})= {\Bbb C}_q[B]$ and
 ${\Bbb C}_q[{\rm{G}}]$}. Let ${\frak g} $ be a semisimple Lie  algebra with  the system of simple roots
 $\alpha_1,\ldots,\alpha_n$. Let ${\rm{G}}$ be its simply connected  Lie group.
  Denote $d_i:=\frac{(\alpha_i,\alpha_i)}{2}$
 and $C:={\Bbb C}[q,q^{-1},(q^{d_i} - q^{-d_i})^{-1}]$.   The quantum universal enveloping
 algebra
is an Hopf algebra over $C$ generated by $E_i, F_i, K_i^{\pm 1}$,
$1\le i\le n$ obeying  Drinfeld-Jimbo relations. The algebra ${\Bbb
C}_q[{\rm{G}}]$,  the subalgebra of the dual Hopf algebra for
$U_q({\frak g})$, is generated by matrix  entries of irreducible
finite dimensional representations $c_{f,v}(a):=f(av)$, $v\in V$,
$f\in V^*$, $a\in U_q({\frak g})$.

We assume that $l$ is admissible. In the case of algebras A3-A4:
 $l$ is admissible if $l$ is odd
and $l\ge 3$ in the case $G$ has ${\rm{G}}_2$ components. The
algebra $C_\varepsilon[{\rm{G}}]$ has a central Poisson subalgebra
$Z_0$ that is isomorphic to ${\Bbb C}[G]$ with  the standard
Belavin-Drinfeld bracket ~\cite{DC-L}. The algebra ${\Bbb
C}_q[{\rm{G}}]$ has a subalgebra $R_q^+$ generated by matrix entries
$c_{f,v_{-\lambda}}$ where $v_{-\lambda}$ is the vector of lowest
weight  in the irreducible representation $V_\lambda$ with highest
weight $\lambda$. The algebra $R_q^+ $ is isomorphic to ${\Bbb
C}_q[{\rm{B}}]$ where ${\rm{B}}:={\rm{B^+}}$. By Drinfeld pairing
the algebra
${\Bbb C}_q[{\rm{B}}]$ is isomorphic to $U_q({\frak b}^-)$. \\
{\bf Proposition 3.4}.  Main theorem is true for
the algebras A3-A4.\\
 {\bf Proof}.
 First, notice that
 the  algebra
${\Bbb C}_q[{\rm{G}}]$ has the denominator subset $S$ generated by
matrix entry $c_{\rho, v_{-\rho}}$ where $\rho$  equal to the half
of sum of positive roots. The localization
 ${\Bbb C}_q[{\rm{G}}]S^{-1}$
is isomorphic to the subalgebra in $U_q({\frak b}^-)\otimes
U_q({\frak b}^+)$ generated by $K_\lambda\otimes K_{-\lambda}$,
$F_i\otimes 1$, $1\otimes E_i$, $1\le i\le n$. It suffices  to
construct an admissible stratification for ${\Bbb C}_q[{\rm{B}}]$.

 The algebra $R_q^+$ (that is equal to $ {\Bbb C}_q[{\rm{B}}]=U_q({\frak b}^-)$) has an  stratification
$({\cal P}_w, S_w)$ where $w$ is an element of the Weyl group $W$
~\cite{DC-L,DC-P2,J}. By definition, the ideal ${\cal P}_w$
 is generated (as ideal) by the elements
 $c_{f,v_{-\lambda}}$ where $f$ is orthogonal to subspace
 $U_q({\frak b}^-)t_wv_{-\lambda}$(here $t_w$  is the corresponding
 element of the braid group). The denominator subset
 $S_w$ is generated  by the element $z_w:=c_{wf^\rho,v_{-\rho}} $
where $f^\rho$ is the element of highest weight in $V_\rho^*$.

 Below we present the other construction of pair $({\cal P}_w, S_w)$.
Decompose the element $w_0$ (of highest length in the $W$) into
product of simple reflections
$$w_0= s_1\ldots s_k s_{k+1}\ldots s_{N}, \quad  s_t:=s_{\alpha_{i_t}}$$
such that $w= s_1\ldots s_k.$ Denote $w_t= s_1\ldots s_t$ (here
$w_k=w$) and $z_t=z_{w_t}$. The elements $z_i$ are  $q$-commute
~\cite[Cor. 3.2]{DC-P2}.  As usual denote $\beta_t:= s_1\ldots
s_{t-1}(\alpha_{i_t})$. The algebra $R_q^+$ is a quantum solvable
algebra with respect to the chain of generators
$$ K^{\pm 1}_1,\ldots, K^{\pm 1}_n,
F_{\beta_1},\ldots,F_{\beta_1},\ldots,F_{\beta_N}.$$

We denote by $B_t$ the subalgebra generated by $ K^{\pm 1}_1,\ldots,
 K^{\pm 1}_n, F_{\beta_k},\ldots,F_{\beta_t}.$ We obtain the
filtration $B_1\subset\ldots\subset B_k\subset B_N=R_q^+$. The
subalgebra $B_k$ depends only on $w$ (denote
$B_w:=B_k$)~\cite{DC-P2}.
 The  element $z_t$ lies in $B_t$ and don't lie in $B_{t-1}$~\cite[Lemma 3.2]{DC-P2}.
Denote by $S_t$ the denominator subset generated by $z_t$.
 The ideal ${\cal P}_w$ has zero intersection with $B_k$ and
  $ B_kS_k^{-1}= R_q^+ S_k^{-1}/ {\cal P}_wS_k^{-1}$.
Let $S_w$ be the denominator subset generated by $S_t$, $1\le t\le
k$ and ${\frak R}_w:= R_q^+S_w^{-1}$. The elements $z_t$ are
FA-elements
 in quantum solvable algebra ${\frak R}_w$~\cite{PanA, PanAA}.
 The adjoint action ${\rm{Ad}}_{z_i}$ are diagonalizable.
 Choose the new generators  $z^{\pm 1}_1,\ldots, z^{\pm 1}_k,
 \tilde{F}_{\beta_{k+1}},\ldots,  \tilde{F}_{\beta_N}$ in localization $R_q^+S_*^{-1}$
 (for definition of $S_* $ see section 2). The ideal
 ${\cal P}_wS_*^{-1}$ is generated by
$ \tilde{F}_{\beta_{t}}$ $k+1\le t\le N$. The elements
$z^l_{t,\varepsilon}$ lie in $Z_0$ ~\cite[Theorem 1.6]{DC-P2}. The
pair $(R_q^+,{\cal P}_w)$ is a standard pair and the stratification
$({\cal P}_w,S_w)$ is an admissible stratification. This verifies
the statement for ${\Bbb C}_q[{\rm{B}}]$ (and therefore for ${\Bbb
C}_q[{\rm{G}}]$).$\Box$

\end{document}